\documentclass[12pt]{article}
\usepackage{amsmath,amsfonts,amssymb,amsthm}

\def\B{\langle}
\def\K{\rangle}

\newcommand{\R}{{\mathbb R}}
\newcommand{\C}{{\mathbb C}}
\newcommand{\Z}{{\mathbb Z}}

\newcommand{\1}{{\bf 1}}

\newcommand{\h}{{\bf h}}
\newcommand{\V}[1]{V_{L}^{#1}}
\newcommand{\Vl}[1]{V_{#1+L}}
\newcommand{\al}{\alpha}
\newcommand{\be}{\beta}
\newcommand{\lm}{\lambda}
\newcommand{\Lm}{\Lambda}
\newcommand{\ga}{\gamma}
\newcommand{\Ce}[1]{{\rm Cent\ }#1}

\newcommand{\eps}{\varepsilon}

\renewcommand{\th}{\theta}

\newcommand{\ci}{\circ}
\newcommand{\ch}{{\rm ch}}

\newcommand{\ot}{\otimes}
\newcommand{\op}{\oplus}
\newcommand{\eqa}{\begin{eqnarray}}
\newcommand{\eeqa}{\end{eqnarray}}
\newcommand{\eqn}{\begin{eqnarray*}}
\newcommand{\eeqn}{\end{eqnarray*}}
\newcommand{\Vn}{V^{\natural}}

\newtheorem{dfn}{Definition}[section]
\newtheorem{pro}[dfn]{Proposition}
\newtheorem{thm}[dfn]{Theorem}

\newtheorem{lem}[dfn]{Lemma}
\newtheorem{cor}[dfn]{Corollary}
\newtheorem{rem}[dfn]{Remark}

\makeatletter
\@addtoreset{equation}{section}

\makeatother

\def\bl{\begin{lem}\label}
\def\el{\end{lem}}
\def\bt{\begin{thm}\label}
\def\et{\end{thm}}
\def\bp{\begin{pro}\label}
\def\ep{\end{pro}}
\def\br{\begin{rem}\label}
\def\er{\end{rem}}
\def\bc{\begin{cor}\label}
\def\ec{\end{cor}}
\def\bd{\begin{dfn}\label}
\def\ed{\end{dfn}}

\def\h{\mathfrak{h}}

\title{\Large Decompositions of the Moonshine Module \\
with respect to subVOAs associated to codes over $\Z_{2k}$}
\author{Hiroki SHIMAKURA}
\date{\small\it Graduate School of Mathematical Sciences,\\
University of Tokyo, Komaba 3-8-1, Tokyo 153-8914, Japan\\
{\rm e-mail: shima@ms.u-tokyo.ac.jp}
}

\begin{document}
\maketitle

\baselineskip 6mm
\begin{abstract}
{In this paper, we give decompositions of the moonshine module $\Vn$ with respect to subVOAs associated to extremal Type II codes over $\Z_{2k}$ for an integer $k\ge2$.
Those subVOAs are isomorphic to the tensor product of 24 copies of the charge conjugation orbifold VOA.
Using such decompositions, we obtain some elements of type $4A$ ($k$ odd) and $2B$ ($k$ even) of the Monster simple group Aut$(\Vn)$.}
\end{abstract}

\section*{Introduction}

The notation of a vertex operator algebra (VOA) is introduced in \cite{Bo1,FLM}.
One of the most interesting examples of a VOA is the moonshine module $\Vn$ constructed in \cite{FLM}.
The automorphism group of $\Vn$ is the Monster, the largest sporadic finite simple group.

The moonshine module $\Vn$ has many subalgebras having good symmetry (\cite{Mi,DLMN}).
The decompositions of $\Vn$ with respect to those subalgebras are computed in \cite{DGH,KLY1,KLY2}.
In particular, $V^{\natural}$ contains a subVOA isomorphic to the tensor product of $48$ copies of $L(1/2,0)$, called a Virasoro frame in \cite{DGH}, and the decomposition of $\Vn$ is computed.

Since $L(1/2,0)$ is isomorphic to $W_2$ algebra at $c=1/2$, we should study the first member of the unitary series of $W_n$ algebras.
For details of $W_n$ algebras, see \cite{FKW} and references given there.
In particular, $W_4$ algebra at $c=1$ is realized as the fixed-point subspace $V_L^+$ of the lattice VOA $V_L$ corresponding to the rank one lattice $L=\Z\al$ with $\B\al,\al\K={6}$ with respect to the $-1$ automorphism of the lattice.
Its fusion rules have a nice symmetry according to \cite{Ab} and for each embedding of $V_L^+$ into $\Vn$, we get a $4A$ element of the Monster \cite{Ma}.

In this paper, we consider more general case $L=\Z\al$ with $\B\al,\al\K=2k$ for an integer $k\ge2$.
Since the Leech lattice $\Lm$ contains many elements of squared length $2k$ and $\Vn$ contains $V_\Lm^+$ as a subVOA, $V^{\natural}$ contains many copies of $V_L^+$.
We consider a set of $24$ mutually orthogonal pairs of opposite vectors of $\Lm$ with squared length $2k$.
We call such a set a $2k$-frame of $\Lm$.
The existence of $2k$-frame of $\Lm$ is shown in \cite{Ch,GH}.
By using a $2k$-frame, we show that $V^{\natural}$ contains a subVOA isomorphic to the tensor product of $24$ copies of $V_L^+$.
Then we give the decomposition of $V^{\natural}$ as a $(V_L^+)^{\ot24}$-module.
By using the fact that there exists the natural bijection between equivalence classes of $2k$-frames and extremal Type II codes of length $24$ over $\Z_{2k}$, the decompositions of $\Vn$ are described in terms of such codes.
By the fusion rules of $V_L^+$, we obtain automorphisms of $\Vn$ with respect to the decompositions.
More precisely, we have $4A$ elements ($k$ is odd) and $2B$ elements ($k$ is even) of the Monster simple group.
Moreover, we give new expressions of McKay-Thompson series for $4A$ elements and obtain formulas of modular functions.

\medskip 
\noindent Professor Ching Hung Lam studies the subject independently.

\medskip

Throughout the paper, we will work over the field $\C$ of complex numbers unless otherwise stated.
We denote the set of integers by $\Z$ and the ring of the integers modulo $k$ by $\Z_k$.

{\it Acknowledgements.} The author wishes to thank Professor Atsushi Matsuo, my research supervisor, for his advice and warm encouragement.
He also thanks Professor Masaaki Kitazume and Professor Ching Hung Lam for helpful advice.

\section{Preliminaries}
In this section, we recall or give some definitions and facts necessary in this paper.

\subsection{Charge conjugation orbifold VOA $\V{+}$ \\ and its irreducible modules}\label{DefofVOA}
In this section, we give the charge conjugation orbifold VOA $\V{+}$ and its irreducible modules.
The description of $\V{+}$ in this paper is slightly different from those given in \cite{DN,Ab} because it is useful to show the main theorem.

Let $L=\Z\al$ be an even lattice of rank one with $\B\al,\al\K=2k$.
Set $H=\C\otimes_{\Z} L$ and we regard $L$ as a subgroup of $H$.
We extend the form $\B\cdot,\cdot\K$ to a $\C$-bilinear form on $H$.
Let $L^{\ci}=\{h\in H\ |\ \B h,\al\K\in\Z\}=L/{2k}$ be the dual lattice of $L$.
Let $\hat{L^\ci}$ be the trivial extension of $L^\ci$ by the order $2$ cyclic group $\B-1\K$.
Form the induced $\hat{L^\ci}$-module $\C\{L^{\ci}\}=\C[\hat{L^\ci}]\ot_{\C[\pm 1]}\C$, where $\C[\cdot]$ denotes the group algebra and $-1$ acts on $\C$ as multiplication by $-1$.
We choose a section $L^\ci\to \hat{L^\ci},\ \be\mapsto e^\be$ such that $e^\be e^\ga=(-1)^{\B\be,\ga\K}e^{\be+\ga}$ for $\be,\ga\in L$.
Set $V_{L^{\ci}}=S(H\otimes t^{-1}\C[t^{-1}])\otimes\C\{L^{\ci}\}$.
For $\lm+L\in L^{\ci}/L$, we set $V_{\lm+L}=S(H\ot t^{-1}\C[t^{-1}])\ot\C\{\lm+L\}$, where $\lm+L$ is regarded as a subset of $L^{\ci}$ and $\C\{\lm+L\}$ is the subspace of $\C\{L^{\ci}\}$ spanned by $\{e^{\mu}\ |\ \mu\in \lm +L\}$.
It is well known that $\V{}$ has a VOA structure and $V_{\lm+L}$ has a $V_L$-module structure for each $\lm+L\in L^{\ci}/L$ (cf. \cite{FLM}).
Note that all the irreducible $\V{}$-modules are given by the set $\{ \Vl{\lm}\ |\ \lm+L \in L^{\ci}/L\}$ (cf. \cite{D1}).
For convenience, we denote $h\ot t^n$ by $h(n)$ for $h\in H$ and $n\in\Z$.

Let $\rho_L$: $L\to\C^{\times}$ be a group homomorphism such that $\rho_L(2\mu)=(-1)^{\B\mu,\mu\K/2}$ for $\mu\in L$.
Let $\th_L$ be the automorphism of $\hat{L}$ given by $\th_L(e^\mu)=\rho_L(2\mu)(e^{\mu})^{-1}$ for $\mu\in L$ and $\th_L(-1)=-1$.
Let $\th_{V_L}$ be the unique commutative algebra automorphism of $V_{L}$ such that $\th_{V_L}(h(n)\ot1)=-h(n)\ot1$ and $\th_{V_L}(1\ot e^{\mu})=1\ot \th_L(e^{\mu})$ for $h\in H$ and $\mu\in L$.

Let $\lm$ be a representative of a coset in $L^{\ci}/L$.
We will extend $\th_{V_L}$ on $V_{L^\ci}$ as a $V_L$-module isomorphism.
Let $\th_{V_L}$: $V_{\lm+L}\to V_{-\lm+L}$ be the $V_L$-module isomorphism such that $\th_{V_L}(e^\lm)=e^{-\lm}$.
When $\lm\in L/2$, the linear map $V_{\lm+L}\to V_{-\lm+L}$, $x\mapsto\rho_L(-2\lm)x$ is an isomorphism of $V_L$-modules.
Then we have the automorphism $\th_{V_L}$ of $V_{\lm+L}$ as a $V_L$-module such that $\th_{V_L}(e^{\lm})=\rho_L(2\lm)e^{-\lm}$.
For each $\th_{V_L}$-stable subspace $M$ of $V_{L^\ci}$, let $M^\pm$ denote the $\pm1$ eigenspaces of $M$ with respect to $\th_{V_L}$ respectively.
Note that $V_L^+$ is a subVOA of $V_L$.

Let $K_L=\{\th_{L}(a)a^{-1}\ |\ a\in\hat{L}\}$ be the subgroup of $\hat{L}$.
Let $T_0$ and $T_1$ be irreducible $\hat{L}/K_L$-modules over $\C$ such that $e^\al$ acts by $\rho_L(\al)$ and $e^\al$ acts by $-\rho_L(\al)$ respectively.
Note that $-1\in\hat{L}$ acts by $-1$ on $T_i$.
Set $\V{T_{i}}=S(H\ot t^{-1/2}\C[t^{-1}])\ot T_{i}$.
Then it has an irreducible $\th_{V_L}$-twisted $\V{}$ module structure.
Moreover, $V_L^{T_0}$ and $V_L^{T_1}$ give all the inequivalent irreducible $\th_{V_L}$-twisted $V_L$-modules (cf. \cite{D2}).

Let $\th_H$ be the unique commutative algebra automorphism of $S(H\ot t^{-1/2}\C[t^{-1}])$ such that $\th_H(h(n))=-h(n)$ for $h\in H$ and $n\in1/2+\Z$.
We consider the linear automorphism of $V_L^{T_i}$ given by $\th_H\ot \1_{T_i}$, where $\1_{T_i}$ is the identity operator of $T_i$.
By abuse of notation, we denote both automorphisms of $V_L^{T_0}$ and $V_L^{T_1}$ by $\th_{V_L^T}$.
Since $e^\al=\th_{V_L}(e^\al)$ on $T_i$, $\th_{V_L^{T}}$ is an automorphism as a $V_L$-module.
We denote the $\pm1$ eigenspaces of $V_L^{T_i}$ with respect to $\th_{V_L^T}$ by $\V{T_i,\pm}$ respectively.
Then $\V{T_{i},\pm}$ becomes an irreducible $\V{+}$-module.

Note that $V_L^+$ and its irreducible modules defined above are isomorphic to those given in {\rm \cite{DN}}.
By Theorem $5.13$ of \cite{DN}, the set
\[
\{\V{\pm},\Vl{\al /2}^{\pm},\V{T_{i},\pm},\Vl{r\al /2k}\ |\ i=0,1,r=1,\dots,k-1\}
\]
gives all inequivalent irreducible $\V{+}$-modules.
\subsection{Moonshine module}\label{DefofMoonshine}
Let us review the moonshine module $V^{\natural}$ from \cite{FLM} for what we need in this paper.

Let $\Lm$ be the Leech lattice with the positive-definite $\Z$-bilinear form $\B\cdotp,\cdotp\K$.
It is a unique positive-definite even unimodular lattice of rank $24$ without roots.

Let $\hat{\Lm}$ be the central extension of $\Lm$ by
the cyclic group $\B -1\K\cong\Z_2$:
\eqa 1\rightarrow\B-1\K\rightarrow\hat{\Lm}~\bar{\to}~\Lm\rightarrow 0\label{cent}
\eeqa
 with the commutator map 
$c_0(\be,\ga)=\B \be,\ga\K+2\Z$ for $\be,\ga\in\Lm$.

Let $\th_{\Lm}$ be the automorphism of $\hat{\Lm}$ given by $\th_{\Lm}(a)=a^{-1}(-1)^{\B \bar{a},\bar{a}\K /2}$ for $a\in\hat{\Lm}$.
We denote the center of $\hat{\Lm}$ by $\Ce{\hat{\Lm}}$.
Set $K_\Lm=\{ \th_{\Lm} (a)a^{-1}\ |\ a\in\hat{\Lm}\}\subset\Ce{\hat{\Lm}}$.
Then $K_\Lm$ is a normal subgroup of $\hat{\Lm}$; consider the group $\hat{\Lm}/K_\Lm$.
Note that the center of ${\hat{\Lm}/K_\Lm}$ is $\B-K_\Lm\K\cong\Z_2$.
By Theorem $5.5.1$ of \cite{FLM}, we have the following proposition.
\bp{T224}$(${\rm \cite{FLM}}$)$\ 
Set $G=\hat{\Lm}/K_\Lm$ and let $\chi$ be a character of $\Ce{G}=\B-K_\Lm\K$ such that $\chi(-K_\Lm)=-1$.
Let $A$ be a maximal abelian subgroup of $G$ and let $\psi$: $A\to\C^{\times}$ be a character with $\psi_{|\Ce{G}}=\chi$.
Then $T=\rm{Ind}_A^G\C_\psi$ is the unique irreducible $G$-module on which $x\in\Ce{G}$ acts by $\chi(x)\1_T$, where $\C_\psi$ is the one-dimensional $A$-module corresponding to $\psi$.
Moreover, $T\cong\oplus\C_{\varphi}$, where $\varphi$ ranges over the characters of $A$ whose restriction to $\Ce{G}$ is $\chi$, and $\dim T=2^{12}$.
\ep

Set the induced $\hat{\Lm}$-module $\C\{\Lm\}=\C[\hat{\Lm}]\ot_{\C[\pm1]}\C$, where $-1$ acts on $\C$ as multiplication by $-1$.
We extend $\th_{\Lm}$ to $\C[\hat{\Lm}]$ linearly.
Since $\th_{\Lm}$ fixes $-1$, we view the automorphism $\th_{\Lm}$ as an automorphism of $\C\{\Lm\}$.

Set $\h=\C\ot_{\Z}\Lm$.
We regard $\Lm$ as a subgroup of $\h$.
We extend $\B\cdot,\cdot\K$ to a $\C$-bilinear form on $\h$.
Let $T$ be the irreducible $\hat{\Lm}/K_\Lm$-module given in Proposition \ref{T224}.
Set $V_{\Lm}=S(\h\ot t^{-1}\C[t^{-1}])\ot\C\{\Lm\}$ and 
$V_{\Lm}^{T}=S(\h\ot t^{-1/2}\C[t^{-1}])\ot T$.
For convenience, we also denote $h\ot t^n$ by $h(n)$ for $h\in\h$ and $n\in\Z/2$.
Let $\th_{V_\Lm}$ be the unique commutative algebra automorphism of $V_{\Lm}$ such that $\th_{V_\Lm}(h(n)\ot 1)=-h(n)\ot 1$ and $\th_{V_\Lm}(1\ot b)=1\ot\th_{\Lm}(b)$, where $h\in\h$, $n\in\Z_{<0}$ and $b\in\C\{\hat{\Lm}\}$.
Let $\th_\h$ be the unique commutative algebra automorphism of $S(\h\ot t^{-1/2}\C[t^{-1}])$ such that $\th_\h(h(n))=-h(n)$ for $h\in\h$, $n\in1/2+\Z_{<0}$.
Then $\th_{V_\Lm^T}=\th_\h\ot (-\1_T)$ is an automorphism of $V_\Lm^{T}$, where $\1_T$ is the identity operator on $T$.
Thus we have $\Vn=V_\Lm^+\op V_\Lm^{T,+}$, where $V_\Lm^+$ is the $\th_{V_\Lm}$-fixed-point subspace of $V_\Lm$ and $V_\Lm^{T,+}$ is the $\th_{V_\Lm^T}$-fixed-point subspace of $V_\Lm^T$.

\subsection{$2k$-frames and codes over $\Z_{2k}$}

In this subsection, we give some terminology on a code over $\Z_{2k}$ (cf. \cite{DHS}).

A ({\it linear}) {\it code} $C$ of length $n$ over $\Z_{2k}$ is a $\Z_{2k}$-submodule of $\Z_{2k}^n$.
We denote the image of $x\in\Z$ with respect to the canonical map $\Z\to\Z_{2k}$ by $\bar{x}$.
We fix an ordered basis of $\Z_{2k}^n$ and denote $i$-th element of this basis by $(\bar0,\dots,\bar0,\bar1,\bar0,\dots,\bar0)$, where $\bar1$ only appears in the $i$-th position.

An element of $C$ is called a {\it codeword}.
The Euclidean weight Ewt$(\cdot)$ on $\Z_{2k}^n$ is given by ${\rm Ewt}(c)=\sum_{i=1}^{24} m_i^2$, where $c=(\bar{m}_1,\dots,\bar{m}_n)\in\Z_{2k}^n$ and $-k<m_i\le k$ for $i=1,\dots,n$.
We define the inner product of $x$ and $y$ in $\Z_{2k}^n$ by $\B x,y\K=\sum_{i=1}^{n}\bar{x}_i \bar{y}_i$, where $x=(\bar{x}_1,\dots,\bar{x}_n)$ and $y=(\bar{y}_1,\dots,\bar{y}_n)$.
The {\it dual code} $C^{\bot}$ of C is defined as $C^{\bot}=\{x\in\Z_{2k}^n\ |\ \B x,y\K =0\ {\rm for\ all}\ y\in C\}$.
$C$ is {\it self-orthogonal} if $C\subset C^{\bot}$ and $C$ is {\it self-dual} if $C=C^{\bot}$.
We define a {\it Type II} code over $\Z_{2k}$ to be a self-dual code with all codewords having Euclidean weight divisible by $4k$.

It is well known that $\tilde{S_n}=\Z_2\wr S_n={\rm Aut}(\Z_{2k}^n)$, where $S_n$ is the symmetric group.
$\tilde{S}_n$ acts on $\Z_{2k}^n$ by the permutation of the coordinate positions and the change of the signs of some positions of $\Z_{2k}^n$.
Two codes $C_0$ and $C_1$ over $\Z_{2k}$ are called {\it equivalent} if they both have length $n$ and if there exists $\sigma\in \tilde{S}_n$ such that $C_0=\sigma(C_1)$.

\bd{D1} {\rm A} {\it {$2k$}-frame} {\rm of a lattice of rank $n$ is a set of $n$ mutually orthogonal pairs of opposite vectors of squared length ${2k}$.}
\ed
We denote the group of isometries of a lattice $\Lm$ by ${\rm Aut}(\Lm)$.
Since ${\rm Aut}(\Lm)$ preserves the inner product, ${\rm Aut}(\Lm)$ acts on the set of ${2k}$-frames of $\Lm$.
We say that ${2k}$-frames $S_0$, $S_1$ of $\Lm$ are {\it equivalent} if there exists $\tau\in{\rm Aut}(\Lm)$ such that $S_0=\tau(S_1)$.

We define an {\it extremal} Type II code of length $24$ over $\Z_{2k}$ to be a Type II code with minimum Euclidean weight $8k$.

\bp{Existex}{\rm \cite{Ch,GH}} For any positive integer $k$, there exists an extremal Type II code of length $24$ over $\Z_{2k}$.
\ep
Using the fact that the Leech lattice $\Lm$ is the unique positive definite unimodular lattice in dimension $24$ without roots, it is easy to see that for $k\ge2$, equivalence classes of $2k$-frames of $\Lm$ are the same as equivalence classes of extremal Type II codes of length $24$ over $\Z_{2k}$.
More precisely, for a $2k$-frame $S$ of $\Lm$, we have an extremal Type II code $C=\Lm/N\subset N^{\ci}/N\cong\Z_{2k}^{24}$, where $N$ is the sublattice of $\Lm$ generated by $S$ and $N^\ci$ is the dual lattice of $N$, and for an extremal Type II code $C$ of length $24$ over $\Z_{2k}$, the lattice constructed by the generalized Construction A with $C$ is the Leech lattice and contains a $2k$-frame.

By direct calculation, we have the following proposition.
\bp{Pcode}
Let $S$ be a ${2k}$-frame of the Leech lattice $\Lm$ and let $N$ be the sublattice of $\Lm$ generated by $S$.
Let $C=\Lm/N$ be a code over $\Z_{2k}$ and let $C_2=(\Lm\cap(N/2))/N$ be a binary code.
Set $m=\dim{C_2}$ over $\Z_2$.
\begin{enumerate}
\item \ $|\Lm/N|={(2k)}^{12}$.
\item \ $m\ge 12$.
\item \ $N/(2\Lm\cap N)$ is an elementary abelian $2$-group with order $2^{24-m}$. 
\item\ If $k$ is odd, then $C_2$ is a Type II code.
\end{enumerate}
\ep

\br{C2code}{\rm Let $S$ be a $2k$-frame of $\Lm$ and $C$ be an extremal Type II code corresponding to $S$.
Let $N$ be the sublattice of $\Lm$ generated by $S$.
Then it is easy to see that $C_2=(\Lm\cap(N/2))/N\cong \{(c_1,\dots,c_{24})\in C\ |\ c_i=0\pmod{k}\ {\rm for\ all\ }i\}$, where we regard the both codes as binary codes.}
\er

\section{Decompositions of $V^{\natural}$ and $\Z_{2\lowercase{k}}$ codes}
In this section, for an integer $k\ge2$, we give the decomposition of $\Vn$ as a $(V_L^+)^{\ot 24}$-module associated with an extremal Type II code over $\Z_{2k}$.
It is easy to see the embedding of $V_L$ into $V_\Lm$.
But the $24$ tensor product of the involution of $V_L$ given in \cite{DN} is not the same involution $\th_{V_\Lm}$ of $V_\Lm$.
By using the definition of $V_L$ given in Section 1.1, we will clarify the problem.

For convenience, we use the following notation.
For $c=(\bar{c}_1,\dots,\bar{c}_{24})\in\Z_{2k}^{24}$, we set
\eqn
M(c)=\bigotimes_{i=1}^{24}V_{c_i\al/{2k}+L}.
\eeqn
Note that $M(c)\cong M(\sigma(c))$ as a $(V_L^+)^{\otimes 24}$-module for $c\in\Z_{2k}^{24}$, where $\sigma$ is an operator which changes signs of some positions.

In particular, for $c=(\bar{c}_1,\dots,\bar{c}_{24})\in\Z_2^{24}$, we set
\eqn
M(c)^+&=& \bigoplus_{e_{i}\in\{\pm\}\atop\prod e_i=+}\bigotimes_{i=1}^{24}\Vl{c_i\al/2}^{e_i},\\
M^T(c)^-&=& \bigoplus_{e_i\in\{\pm\}\atop\prod e_i=-}\bigotimes_{i=1}^{24} V_{L}^{T_{c_i},e_i,},
\eeqn
where $\{\pm\}\cong\B -\K\cong\Z_2$.

\br{R30} {\rm Let $c=(c_1,\dots,c_{24})$ be an element of $\Z_2^{24}$.
\begin{enumerate}
\item $M(c)^+$ is the fixed-point subspace of $\ot V_{c_i\al/2+L}$ with respect to $\th_{V_L}^{\ot24}$.
\item $M^T(c)^-$ is the $-1$ eigenspace of $\ot V_L^{T_{c_i}}$ with respect to $\th_{V_L^T}^{\ot24}$.
\end{enumerate}}
\er

\subsection{Main results}
The following is our main theorem.
\bt{T2} For $k\ge2$, let $C$ be an extremal Type II code of length $24$ over $\Z_{2k}$ and set $C_2=\{(c_1,\dots,c_{24})\in C\ |\ c_i=0\pmod{k}\ {\rm for\ all\ }i\}$, which is binary code. 
Set $m=\dim C_2$ over $\Z_2$.
Then there exists an embedding of $(\V{+})^{\ot24}$ into $V_{\Lm}^{+}$ as a subVOA such that $V^{\natural}=V_{\Lm}^+\op V_{\Lm}^{T,+}$ decomposes into $(\V{+})^{\ot24}$-modules as the following:
\begin{enumerate}
\item \eqn
V_{\Lm}^+\cong\bigoplus_{c\in C_2} M(c)^+\op\frac{1}{2}\bigoplus_{c\in C\setminus C_2}M({c}).
\eeqn
\item \eqn
V_{\Lm}^{T,+}\cong\bigoplus_{c\in C_2^\bot}2^{m-12}M^T(c)^-.
\eeqn
\end{enumerate}
\et
\br{R321}{\rm This decomposition is uniquely determined by the extremal Type II code of length $24$ over $\Z_{2k}$, up to the action of $\tilde{S_n}=\Z_2\wr S_n$.}
\er
The rest of this section is devoted to the proof of the theorem.

\subsection{Key lemma}

In this subsection, we will give the key lemma for Theorem \ref{T2}.
It is an extension of Theorem D.6 of \cite{DGH}.
More precisely, we consider the case of a dual lattice.
We will use the lemma to identify $\th_{V_L}^{\ot24}$ with $\th_{V_\Lm}$ on $V_\Lm$.

Let $U$ be an even integral lattice and let $V_U$ be the lattice VOA associated to $U$.
By \cite{D1}, $V_{U^\ci}$ is a $V_U$-module.
We choose a section $U^\ci\to\hat{U^\ci}$, $x\mapsto e^x$ such that the $2$-cocycle with respect to it is $\Z$-bilinear.
\bd{lift} A {\it lift of} $-1$ of $U^\ci$ is an automorphism $\th$ of $V_{U^\ci}$ as a $V_U$-module such that for all $x\in U^\ci$, there is a scalar $c_x$ so that $\th:e^x\mapsto c_xe^{-x}$.
\ed
Set $A_U={\rm Hom}(U,\C^\times)$ and set $A_{U^\ci}=\{f:U^\ci\to\C^\times\ |\ f(x+y)=f(x)f(y)\ {\rm  for}\ x\in U,\ y\in U^\ci ,\ f_{|U}\in A_U\}$.
For $f\in A_{U^\ci}$, we set the automorphism $\tilde{f}$ of $V_{U^\ci}$ as a $V_U$-module by $\tilde{f}$: $y\ot e^x\mapsto f(x)y\ot e^x$, and we set $\tilde{A}_{U^\ci}=\{\tilde{f}\ |\ f\in A_{U^\ci}\}$.

By Theorem D.6 of \cite{DGH}, we have the following proposition.
\bp{conjuoflift1}For lifts $f$ and $g$ of $-1$ of $U^\ci$, there exists an element $s\in\tilde{A}_{U^\ci}$ such that $s^{-1}fs_{|V_U}=g_{|V_U}$.
\ep

Moreover, we consider the case of $V_{U^\ci}$.
\bl{conjuoflift2}
Let $f,g$ be lifts of $-1$ of $U^\ci$ such that $f=g$ on $V_U$.
Let $S$ be a set of representatives of $(U^\ci\cap U/2)/U$.
Then $f$ and $g$ are conjugate by an element of $\tilde{A}_{U^\ci}$ whose restriction on $V_U$ is the identity map if and only if $f(e^x)=g(e^x)$ for $x\in S$.
\el
\begin{proof}\ In order to simplify the proof, we assume $g(e^x)=e^{-x}$ for $x\in U^\ci$.
Let $\tilde{S}$ be a set of representatives of $U^\ci/U$ containing $S$.

First, we assume $f(e^x)=g(e^x)$ for $x\in \tilde{S}$.
By a conjugation of $\tilde{A}_{U^\ci}$, we assume $f(e^x)=e^{-x}$ for $x\in\tilde{S}$.
Namely, if $f(e^x)=c_xe^{-x}$ for $x\in S$, we set $s\in A_{U^\ci}$ such that $s(x)=c_x^{-1/2}$ for $x\in \tilde{S}\setminus S$ and $s(y)=1$ for $y\in U\cup S$.
It is easy to see that $\tilde{s}^{-1}f\tilde{s}(e^x)=e^{-x}$ for $x\in \tilde{S}$.
Note that $V_{U^\ci}=\op_{\lm\in\tilde{S}} V_{\lm+U}$, and for $\lm\in\tilde{S}$, $V_{\lm+U}$ is generated by $e^{\lm}$ as a $V_U$-module.
Since $f$ is $V_U$-module automorphism and $2$-cocycle is a $\Z$-bilinear map, we have $f(e^x)=e^{-x}$ for $x\in U^\ci$ by the direct calculation.

Next, we assume $\tilde{s}^{-1}f\tilde{s}=g$ for $\tilde{s}\in\tilde{A}_{U^\ci}$ such that $\tilde{s}_{|V_U}$ is the identity map.
For $x\in S$, we set $f(e^x)=c_xe^{-x}$, where $c_x\in\C^\times$.
Then we have $\tilde{s}^{-1}f\tilde{s}(e^x)=s(x)c_x\tilde{s}^{-1}(e^{-2x+x})=s(-2x)^{-1}c_xe^{-x}=s(-2x)g(e^x)$ for $x\in S$.
Since $-2x\in U$ and $s=1$ on $U$, we have $f(e^x)=g(e^x)$ for $x\in \tilde{S}$.
\end{proof}

\subsection{Decomposition of the untwisted space $V_{\Lm}^+$}\label{S3.1}
\noindent In this subsection, we give the decomposition of $V_{\Lm}^+$.

Let $F=(x_1,\dots,x_{24})$ be the $2k$-frame of $\Lm$ corresponding to the code $C$ and let $N$ be the sublattice of $\Lm$ generated by $F$.
Since $L=\Z\alpha$ is an even lattice of rank one with $\B\alpha,\alpha\K=2k$, we have $N\cong L^{\op24}$.
From \cite{D1}, it follows that $V_{\Lm}$ is decomposed as
\eqa
V_{\Lm}=\bigoplus_{\lm+N\in\Lm /N}V_{\lm +N}\label{DOU}
\eeqa
as a $V_N$-module.

It is well known that
\eqa
V_{\lm+N}\cong\bigotimes_{i=1}^{24}\Vl{c_i\al/{2k}}= M(c)\label{dop}
\eeqa
as a $(V_L^+)^{\ot 24}$-module, where $\lm +N =\sum_{i=1}^{24}c_ix_i/{2k}+N$ and $c=(\bar{c}_1,\dots,\bar{c}_{24})$ is a codeword of $C=\Lm/N$.

Now, we have two involutions of $V_\Lm$ as a $V_N$-module.
One is $f=\th_{V_\Lm}$ given in Section \ref{DefofMoonshine} and the other one is $g=\th_{V_L}^{\ot24}$.
Note that $\th_{V_L}$ is the automorphism of $V_{L^\ci}$ as a $V_L$-module given in Section \ref{DefofVOA}.
We will apply Proposition \ref{conjuoflift1} and Lemma \ref{conjuoflift2} to $V_\Lm\subset V_{N^\ci}$, and identify $f$ and $g$.
Since we consider $V_\Lm$ in this section, we use the notation $\tilde{A}_{\Lm}$ instead of $\tilde{A}_{N^\ci}$.
By \cite{FLM}, we can choose a $\Z$-bilinear $2$-cocycle $\eps:\Lm\times\Lm\to\Z_2$ such that $\eps(x,x)=\B x,x\K/2$.

By Proposition \ref{conjuoflift1} there exists $\tilde{s}\in \tilde{A}_{\Lm}$ such that $\tilde{s}^{-1}f\tilde{s}=g$ on $V_N$.
Therefore we identify $g$ with $f$ on $V_N$, and we have the inclusion of subVOAs $(V_L^+)^{\ot24}\subset V_\Lm^+\subset V^{\natural}$.

\br{R3} {\rm Let $V_{L,\R}^+$ and $\Vn_\R$ be the real forms of $V_L^+$ and $\Vn$ as constructed in \cite{FLM}.
By \rm{\cite{FLM}}, those have the positive-definite symmetric bilinear forms.
By the above inclusion, we have the embedding $V_{L,\R}^+\subset \Vn_\R$.
In this embedding, the positive-definite symmetric bilinear form of the subVOA is the restriction of that of $\Vn_\R$.}
\er
Next, in order to apply to lemma \ref{conjuoflift2}, we have to check the hypothesis of the lemma.
Let $s$ be the element of $A_\Lm$ corresponding to $\tilde{s}$.
For $x\in \Lm\cap N/2$, we have $\tilde{s}^{-1}f\tilde{s}(e^x)=s(-2x)e^{-x}$ and $g(e^x)=\rho_N(2x)e^{-x}$, where $\rho_N=\op\rho_{\Z x_i}$.
Since $\rho_N$ and $s$ are linear on $N$, we can choose $\rho_N$ and the set $S$ of representatives of $(\Lm\cap N/2)/N$ such that $\tilde{s}^{-1}f\tilde{s}=g$ on $\tilde{S}$.
Therefore we can use lemma \ref{conjuoflift2}, and we have $f$ and $g$ are conjugate by $\tilde{A}_{\Lm}$.
So, we identify $f$ with $g$ on $V_\Lm$.

We consider the $\th_{V_L}^{\ot24}$-fixed-point subspace of $V_\Lm$.
We obtain the following lemma (cf. \cite{KLY1}).
\bl{L2}Let $c$ be a codeword of $C=\Lm/N$.
\begin{enumerate}
\item If $c\in C_2=(\Lm\cap N/2)/N$, then $M(c)$ is $\th_{V_L}^{\ot24}$-invariant, and $M(c)^{+}$ is the $\th_{V_L}^{\ot24}$-fixed-point subspace.\\
\item If $c\in C\setminus C_2$, then $M(c)\op M(-c)$ is $\th_{V_L}^{\ot24}$-invariant, and the $\th_{V_L}^{\ot24}$-fixed-point subspace 
$(M(c)\op M(-c))^{+}$ is isomorphic to $M(c)$ as a $(V_L^+)^{\ot24}$-module.
\end{enumerate}
\el

By Lemma \ref{L2}, (\ref{DOU}) and (\ref{dop}), we obtain Theorem \ref{T2} (i).
\subsection{Decomposition of the twisted space $V_{\Lambda}^{T,+}$}
In this subsection, we give the decomposition of $V_{\Lm}^{T,+}$.

Set $K(N)=K_\Lm\cap\hat{N}$ and $K(2N)=K_\Lm\cap\hat{2N}$.
Since $T$ is a $\hat{\Lm}/K_\Lm$-module, $T$ is a $\hat{N}K_\Lm/K_\Lm$-module.
Note that $-1\in\hat{\Lm}$ acts on $T$ as multiplication by $-1$.
Let $M$ be a maximal abelian subgroup of $\hat{\Lm}/K_\Lm$ such that $M\supset \hat{N}K_\Lm/K_\Lm$.
Note that $|M|=2^{13}$.
By Proposition \ref{T224}, $T$ decomposes into the direct sum of all the irreducible $M$-modules on which $-K_\Lm$ acts by $-1$.
Therefore $T$ decomposes into the direct sum of irreducible $\hat{N}K_\Lm/K_\Lm$-modules on which $-K_\Lm$ acts by $-1$.
By Proposition \ref{Pcode} ($3$), we have $|\hat{N}K_\Lm/K_\Lm|=2|N/(N\cap2\Lm)|=2^{25-m}$.
Then any multiplicities of the irreducible $\hat{N}K_\Lm/K_\Lm$-module is $2^{m-12}$.
Since $\hat{N}K_\Lm/K_\Lm\cong\hat{N}/K(N)$, $T$ is a $\hat{N}/K(N)$-module.
Thus we obtain
\eqa
T\cong\bigoplus_{\psi\in {\rm Irr}(\hat{N}/K(N))\atop\psi(-K(N))=-1}2^{m-12}T_{\psi}, \label{deo1}
\eeqa
where Irr$(\hat{N}/K(N))$ is the set of the characters for $\hat{N}/K(N)$ and $T_{\psi}$ is the one-dimensional module $\C$ with a character $\psi$.

By using the canonical map $\pi$: $\hat{N}/K(2N)\to \hat{N}/K({N})$, we view $T$ as a $\hat{N}/K(2N)$-module.
Note that the kernel of $\pi$ is $K(N)/K(2N)$.
Therefore, as a $\hat{N}/K(2N)$-module, $T$ is decomposed
\eqa
T\cong\bigoplus_{\chi\in {\rm Irr}(\hat{N}/K(2N))\atop\chi(-K(2N))=-1}m_{\chi}T_{\chi},\label{dot}
\eeqa
where $m_{\chi}$ is the multiplicity.
By ($\ref{deo1}$), we have $m_{\chi}\in\{0,2^{m-12}\}$.
By (\ref{dot}), $V_{\Lm}^{T}$ has a decomposition,
\eqa
V_{\Lm}^{T}\cong \bigoplus_{\chi\in {\rm Irr}(\hat{N}/K(2N))\atop\chi(-K(2N))=-1}m_{\chi}V_{N}^{T_{\chi}},\label{dt1}
\eeqa
where $V_{N}^{T_{\chi}}$ is the $\th_{V_\Lm}$-twisted $V_{N}$-module corresponding to the $\hat{N}/K(2N)$-module $T_{\chi}$.

As in Section 2.3, we fix the ordered basis of $N$ as $(x_1,\dots,x_{24})$.
Let $L_i=\Z x_i$ be a sublattice of $N$ and set $K(L_i)=K_\Lm\cap \hat{L_i}$.
For $\chi\in{\rm Irr}(\hat{N}/K(2N))$, we set the character of $\hat{L_i}/K(L_i)$ $\chi_i(e^{x_{i}})=\chi(e^{x_{i}})$.
For a character $\chi=\ot_{i=1}^{24}\chi_{i}$ such that $\chi(-K(2N))=-1$, we set $P(\chi)=(c_1,\dots,c_{24})\in\Z_2^{24}$, where $\chi_i(\tilde{s}(e^{x_i}))=\rho_i(x_i)(-1)^{c_i}$.
For $x=\sum d_ix_i\in N$, we set $Q(x)=(\bar{d}_1,\dots,\bar{d}_{24})\in\Z_2^{24}$.

Then we have 
\eqn
\chi(\tilde{s}(e^x))=\rho(\sum d_ix_i)(-1)^{\sum c_id_i}=\rho(x)(-1)^{\B P(\chi),Q(x)\K}.
\eeqn

Now, let us determine the multiplicities $\{m_{\chi}\}$.
Let $\chi$ be an element of ${\rm Irr}(\hat{N}/K(2N))$ such that $\chi(-K(2N))=-1$.
Set $n=2^{m-12}$.
Then $m_{\chi}=n$ if and only if there exists $\tilde{\chi}\in{\rm Irr}(\hat{N}/K(N))$ such that $\chi=\tilde{\chi}\ci\pi$, namely, $\chi(x)=1$ for any $x\in {\rm Ker}\pi=K(N)/K(2N)$.
Since $\th_\Lm(e^x)(e^x)^{-1}=e^{-x}(-1)^{\eps(x,-x)}e^{-x}=e^{-2x}=\rho_N(2x)s(-2x)e^{-2x}$ for $x\in \Lm\cap N/2$, we have $K(N)=\{\th_{V_\Lm}(a)a^{-1}\ |\ \bar{a}\in\Lm\cap N/2\}=\{\rho_N(-\be)\tilde{s}(e^{\be})\ |\ \be\in2\Lm\cap N\}$.
Therefore $\chi(x)=1$ for any $x\in K(N)/K(2N)$ if and only if $\B P(\chi),Q(x)\K=0$ for any $x\in K(N)/K(2N)$

Since $K(N)/K(2N)\cong(N\cap2\Lm)/2N\cong (\Lm\cap (N/2))/N=C_2$, $m_{\chi}=n$ if and only if $P(\chi)\in C_2^{\bot}$.
By (\ref{dt1}), we have an isomorphism 
\eqa
V_\Lm^T\cong2^{m-12}\bigoplus_{\chi\in{\rm Irr}(\hat{N}/K(2N))\atop {\chi(-K(2N))=-1\atop P(\chi)\in C_2^\bot}}V_N^{T_{\chi}}\label{mult}.
\eeqa

Next, we fix a character $\chi$ and consider the space $V_{N}^{T_{\chi}}$.
If $P(\chi)=(c_1,\dots,c_{24})$, then $T_{\chi}=\ot T_{c_i}$.
Thus, we get an isomorphism $V_{N}^{T_{\chi}}\cong\ot V_L^{T_{c_i}}$.
Note that we have $\th_{V_\Lm^T} = -(\th_{V_L^T}^{\ot24})$, because $\th_{V_\Lm^T}$ acts by $-1$ on $T_{\chi}$ and $\th_{V_L^{T}}$ acts by $1$ on $T_{c_{i}}$.
Therefore, for $\chi\in {\rm Irr}(\hat{N}/K(2N))$ with $\chi(-K(2N))=-1$, we have an isomorphism of 
$(V_{L}^{+})^{\ot 24}$-modules
\eqa
V_{N}^{T_{\chi},+}\cong M^T(P(\chi))^{-}.
\label{dt2}
\eeqa
By $(\ref{mult})$ and $(\ref{dt2})$, we get Theorem \ref{T2} (ii).

\section{Character and automorphisms of the moonshine module}
In this section, we give the character of $\Vn$ and give some automorphisms of $\Vn$.
\subsection{Character of the moonshine module}

We recall the character of a VOA.
Let $V=\op_{n=0}^{\infty}V_n$ be a VOA and the character of $V$ is given by $\ch(V)=\sum_{n=0}^{\infty}\dim (V_n) q^n$.

In order to give the character of $\Vn$, we consider the symmetrized weight enumerator of a code over $\Z_{2k}$.
The symmetrized weight enumerator of a code $C$ over $\Z_{2k}$ is defined as:
\eqn
{\rm swe}_C(p_0,\dots,p_k)=\sum_{c\in C}p_0^{n_0(c)}p_1^{n_1(c)}\dots p_k^{n_k(c)},
\eeqn
where $n_i(c)$ denotes the number of $j$ such that $c_j=\pm i$.
Note that the symmetrized weight enumerators of equivalence codes are same.
For $0\le i\le k$, we set
\eqn
 a_i&=&{\rm ch} V_{i\al/2k+L}= \frac{1}{\phi(q)}\sum_{j\in \Z}q^{k(j+i/2k)^2}, \\
  b&=&({\rm ch}V_{L}^{+}-{\rm ch}V_{L}^{-})=\frac{1}{\phi(q)}\sum_{j\in\Z}(-1)^jq^{j^2},
\eeqn
where $\phi(q)=\prod_{n\ge1}(1-q^n)$.
Note that for $c\in C_2\setminus\{0\}$, we have $\ch{(M(c)^+)}=\ch{(M(c))}$/2, and $\ch{(M(0)^+)}=b^{24}+\ch{(M(0))}/2$. 

By Theorem \ref{T2}, we obtain the following corollary.
\bc{Jfunc} Let $C$ be an extremal Type {\rm II} code of length $24$ over $\Z_{2k}$.
Then we have
\eqn
{\rm ch}(\Vn)=q(J(q)-744)&=&\frac{1}{2}{\rm swe}_C(a_0,a_1,\dots,a_k)+\frac{1}{2}b^{24}\\ &+&2^{11}q^{3/2}\big(\frac{\phi(q)^{24}}{\phi(q^{1/2})^{24}}-
\frac{\phi(q^2)^{24}\phi(q^{1/2})^{24}}{\phi(q)^{48}}\big).
\eeqn
\ec

\br{ChFLM}{\rm It is easy to see that ${\rm swe}_C(a_0,a_1,\dots,a_k)=\Theta_\Lm(q)/(\phi(q))^{24}$, where $\Theta_\Lm(q)$ is the theta function associated with $\Lm$.
Since $b=\phi(q)/\phi(q^2)$, we obtain}
\eqn
{\rm ch}(\Vn)=\frac{1}{2}\big(\frac{\Theta_\Lm(q)}{\phi(q)^{24}}+\frac{\phi(q)^{24}}{\phi(q^2)^{24}}\big)+2^{11}q^{3/2}\big(\frac{\phi(q)^{24}}{\phi(q^{1/2})^{24}}-\frac{\phi(q^2)^{24}\phi(q^{1/2})^{24}}{\phi(q)^{48}}\big).
\eeqn
{\rm This equation is given in Remark $10.5.8$ of \rm{\cite{FLM}}.}
\er

\subsection{Automorphisms of the moonshine module}\label{AMM}

Let $\Gamma$ be the set of inequivalent irreducible $V_L^+$-modules.
If $k$ is odd then we define the map $\mu_k:\Gamma\to\C^\times$ setting by
\[
\mu_k(W)=
\left\{\begin{array}{cl}
 \mbox{$1$} & \mbox{${\rm if}\ W\in{\{\V{\pm}, \Vl{j\al /k}\ |\ 1\le j\le (k-1)/2\}}$},\\
 \mbox{$-1$} & \mbox{${\rm if}\ W\in \{\Vl{\al /2}^{\pm}, \Vl{(2j-1)\al /2k}\ |\ 1\le j\le (k-1)/2\}$},\\
 \mbox{$i$} & \mbox{${\rm if}\ W\in\ \{\V{T_{0},\pm}\}$},\\
 \mbox{$-i$} & \mbox{${\rm if}\ W\in\ \{\V{T_{1},\pm}\}$},
\end{array}
\right.
\]
and if $k$ is even then we defined the map $\mu_k:\Gamma\to\C^\times$ setting by
\[
\mu_k(W)=
\left\{\begin{array}{cl}
 \mbox{$1$} & \mbox{${\rm if}\ W\in\{\V{\pm}, \Vl{\al /2}^{\pm},\Vl{j\al /2k}\ |\ 1\le j\le (k-1)\}$},\\
 \mbox{$-1$} & \mbox{${\rm if}\ W\in \{\V{T_{0},\pm}, \V{T_{1},\pm}\}$}.
\end{array}
\right.
\]

The fusion algebra of $\V{+}$ is the vector space $U=\op_{W\in\Gamma}\C W$ equipped products $\times$ given by fusion rules, where we regard $W$ as a formal element.
The definition of fusion rules is given in \cite{DL}.
An automorphism of the fusion algebra $U$ is a linear automorphism $g$ such that $g(A\times B)=g(A)\times g(B)$ for $A,B\in U$.
By the fusion rules of $\V{+}$ determined in \cite{Ab}, we have the following proposition (cf. \cite{Ma}).
\bp{AOF}
The linear map of the fusion algebra of $\V{+}$, $W\mapsto\mu_k(W)W$, is an automorphism of the fusion algebra of $\V{+}$. 
\ep
Suppose we are given a decomposition 
\eqn
V^{\natural}=\bigoplus_{W_{i}\in \Gamma}m_{W_1,\dots,W_{24}}W_1\ot\cdots\ot W_{24}
\eeqn
as a $(V_L^+)^{\ot 24}$-module, where $m_{W_1,\dots,W_{24}}$ is the multiplicity.
For each $i\in\{1,\dots,24\}$, we define a linear automorphism $\sigma_i$ of $V^\natural$ by
\eqn
\sigma_{i}(x)=\mu_k(W_i)x
\eeqn
for $x\in\bigotimes_{j=1}^{24} W_j$.
By Proposition \ref{AOF}, we have the following proposition (cf. \cite{Mi}).
\bp{4A}
$\sigma_{i}$ is a VOA automorphism of $V^{\natural}$.
\ep
\begin{proof}\ By the definition of $\mu_k$, $\sigma_i$ fixes any element of $(V_L^+)^{\ot 24}$.
In particular, $\sigma_i$ fixes the Virasoro element and the vacuum vector.
Note that the fusion rules of tensor products of modules are the tensor products of the fusion rules of those modules.

Let $\Gamma_{24}=\{\ot_{i=1}^{24}W_i\ |\ W_i\in\Gamma\}$.
Let $\mu_{k,i}$: $\Gamma_{24}\to\C^\times$ be a map such that $\mu_{k,i}(W)=\mu_k(W_i)$ for $W=\ot W_j\in\Gamma_{24}$.
Let $W^1$, $W^2$ be elements of $\Gamma_{24}$.
By the definition of fusion rules, we have $Y(v,z)w\in (W^1\times W^2) [z,z^{-1}]$ for $v\in W^1$, $w\in W^2$, where $Y(\cdot,z)$ is the vertex operator of $V^{\natural}$ and $\times$ is fusion rules for $V_L^+$.
Therefore we have $Y(\sigma_i(v),z)\sigma_i(w)=\mu_{k,i}(W^1)\mu_{k,i}(W^2)Y(v,z)w=\mu_{k,i}(W^1\times W^2)Y(v,z)w=\sigma_i(Y(v,z)w)$.
\end{proof}
By Theorem \ref{T2}, we have the following proposition.
\bp{RRR4} Suppose $k$ is odd.
In decompositions of $\Vn$ given by Theorem $\ref{T2}$, $\sigma_{i}$ is a $4A$ element of the Monster.
In fact, $\sigma_i\in 2_{+}^{1+24}\subset 2_+^{1+24}.{\it Conway_{1}}\subset {\it Monster}$, where $2_+^{1+24}.{\it Conway_{1}}$ is a non-split extension of ${\it Conway_1}$ $($largest simple Conway group$)$ by the extra-special group $2_+^{1+24}$.
\ep
\begin{proof}\ By Theorem \ref{T2} and Proposition \ref{AOF}, $\sigma_i^2$ acts by $1$ on $V_\Lm^+$ and acts by $-1$ on $V_\Lm^{T,+}$.
By \cite{FLM}, the centralizer of $\sigma_i^2$ in the Monster simple group is $2_+^{1+24}.{\it Conway_{1}}$.
Since $\sigma_i$ commutes with $\sigma_i^2$, we have $\sigma_i\in 2_+^{1+24}.{\it Conway_{1}}$.
Since $\sigma_i$ preserves the space $\C e^{\be}$ for any $\be\in\hat{\Lm}$, we have $\sigma_i\in 2_+^{1+24}$.
It is well known that there are the four types of Monster elements $1A$, $2A$, $2B$ and $4A$ contained in $2_+^{1+24}$.
Since $\sigma_i$ is an order $4$ element, $\sigma_i$ is a $4A$ element of the Monster.\end{proof}

\br{2B}{\rm Suppose $k$ is even.
In decompositions of $\Vn$ given by Theorem \ref{T2}, $\sigma_i$ is a $2B$ element of the Monster.
More precisely, $\sigma_i$ acts by $1$ on $V_{\Lm}$ and acts by $-1$ on $V_{\Lm}^T$, namely, $\sigma_i=z$ given in $(10.4.48)$ of {\rm \cite{FLM}}.}
\er

\subsection{McKay-Thompson series for 4A elements of the Monster}
In this section, we assume that $k$ is odd.
Then we give the McKay-Thompson series for 4A elements $\sigma_i$ given in Section 3.2.
The expressions of it are different from \cite{CN}, and we obtain formulas of modular functions.

We recall the McKay-Thompson series.
Let $g$ be an element of the Monster simple group and let $\Vn=\sum_{n=0}^{\infty}\Vn_n$ be the moonshine module.
Then the McKay-Thompson series for $g$ is given by $T_{g}(q)=q^{-1}(\sum_{n=0}^{\infty}({\rm Tr}\ g_{|\Vn_n})q^n)$, where ${\rm Tr}\ g_{|\Vn_n}$ is the trace of the action of $g$ on $V_n$.

Let $\tau_i$ :$\Z_{2k}^{24}\to\Z_{2}$ be the composite map of the projection map $\Z_{2k}^{24}\to\Z_{2k}$ with respect to $i$-th elements and the canonical map $\Z_{2k}\to\Z_{2}$.
Note that for $c\in C$, the automorphism $\sigma_i$ acts by the multiplication $(-1)^{\tau_i(c)}$ on $M(c)$.

\bl{MT} Let $\sigma_i$ be the automorphism of $\Vn$ given in Section 3.2.
In the decomposition given in Theorem 2.1, we have
\begin{enumerate}
\item \eqn
\sum_{n=0}^{\infty}({\rm Tr}\ {\sigma_i}_{|(V_\Lm^+)_n})q^n=\frac{1}{2}\sum_{c\in C}(-1)^{\tau_i(c)}\ch{(M(c))}+b^{24},
\eeqn
where $b$ is given in Section 3.1.
\item \eqn
\sum_{n=0}^{\infty}({\rm Tr}\ {\sigma_i}_{|(V_\Lm^{T,+})_n})q^n&=&0.
\eeqn
\end{enumerate}
\el
\begin{proof} By direct calculation, we have (i).
Since $C_2$ contains the all $1$ element, the numbers of elements of $C_2$ whose $i$-th coordinate is $0$, and whose $i$-th coordinate is $1$ are equal.
By the definitions of $V_L^{T_i}$, we have $\ch{(V_L^{T_0,\pm})}=\ch{(V_L^{T_1,\pm})}$ respectively.
Therefore we have (ii).
\end{proof}

\bc{McK}
Let $C$ be an extremal Type II code over $\Z_{2k}$.
The McKay-Thompson series for the $4A$ element $\sigma_i$ is given by 
\eqn
T_{4A}(q)=q^{-1}\{\frac{1}{2}\sum_{c\in C}(-1)^{\tau_i(c)}\ch{(M(c))}+b^{24}\}.
\eeqn
\ec

\br{Modular}{\rm In \cite{Bo2,CN}, all the McKay-Thompson series are computed.
So we have the following formulas of modular functions.
\eqn
\frac{\eta(q^2)^{48}}{\eta(q)^{24}\eta(q^4)^{24}}-24=q^{-1}\frac{1}{2}\sum_{c\in C}(-1)^{\tau_i(c)}\ch{(M(c))}+\frac{\eta(q)^{24}}{\eta(q^2)^{24}},
\eeqn
where the Dedekind $\eta$-function $\eta(q)=q^{\frac{1}{24}}\prod_{n\ge1}(1-q^n)$.}
\er

\end{document}